\documentclass[11pt]{article}

\usepackage{fullpage,amsmath,amssymb,color}
\usepackage{algorithm}
\usepackage{enumerate}
\usepackage{algorithmic}
\newcommand{\BlackBox}{\rule{1.5ex}{1.5ex}}

\usepackage{graphicx}
\usepackage{epstopdf}
\usepackage{subfigure}
\graphicspath{{figures/}}

\newcommand{\sA}{{\sigma_{\min}(A)}}
\newcommand{\sB}{{\sigma_{\max}(B)}}
\newcommand{\sG}{{\sigma_{\max}(G)}}
\newcommand{\smH}{{\sigma_{\min}(H)}}
\newcommand{\sMH}{{\sigma_{\max}({H})}}
\newcommand{\stH}{{\sigma_{\max}(\tilde{H})}}
\newcommand{\cG}{{\cal G}}

\newtheorem{corollary}{Corollary}[section]
\newtheorem{theorem}{Theorem}[section]

\newcommand{\Real}{{\mathbb R}}
\newcommand{\cL}{{\mathcal L}}

\title{On the Duality Gap Convergence of ADMM Methods}
\date{}
\author{
Da Tang \\
Tsinghua University \\ 
 Beijing, China
\and 
Tong Zhang \\
Baidu Inc, Beijing, China \\
Rutgers University, NJ, USA
}

\begin{document}
\maketitle

\begin{abstract}
This paper provides a duality gap convergence analysis for the standard ADMM
as well as a
linearized version of ADMM. 
It is shown that under appropriate
conditions, both methods achieve linear
convergence. However, the standard ADMM achieves a
faster accelerated convergence rate than that of the linearized ADMM. 
A simple numerical example is used to illustrate the difference in convergence behavior.
\end{abstract}

\section{Introduction}
This paper considers the following optimization problem:
\begin{equation}
\begin{aligned}
\min\limits_{w,v}&\ [ \phi(w)+g(v)] \\
\text{subject to}\ &Aw-Bv=c ,
\end{aligned}  
\label{eqn:aecco}
\end{equation}
where $(w,v)\in \Real^n \times \Real^m$ are unknown vectors, $A\in \Real^{p\times
  n}$, $B\in \Real^{p\times m}$ and $c\in\Real^p$ are known matrices and
vector. 
In this paper, we assume that $\phi:\mathbb R^n\rightarrow \mathbb
R\cup\{+\infty\}$ and $g:\mathbb R^m\rightarrow \mathbb
R\cup\{+\infty\}$ are convex functions. 

A popular method for solving \eqref{eqn:aecco} is the Alternating
Direction Method of Multipliers (ADMM) algorithm. It solves the
problem by alternatively optimizing the variables in the Augmented Lagrangian function:
\begin{eqnarray}
\cL(w,v,\alpha,\rho)=\phi(w)+g(v) + \alpha^\top
(Aw-Bv-c)+\frac{\rho}{2}\|Aw-Bv-c\|_2^2 , \label{eq:lagrangian}
\end{eqnarray}
and the resulting procedure is summarized in
Algorithm~\ref{algo:admm}. In the algorithm, both
$G$ and $H$ are symmetric positive semi-definite matrices.
In the standard ADMM, we can set $G=0$ and $H=0$. 
The method of introducing the additional term 
$\|v-v^{t-1}\|_G^2 = (v-v^{t-1})^\top G (v-v^{t-1})$ 
is often referred to as preconditioning. If we let 
$G= \beta I - B^\top B$ 
for a sufficiently large $\beta>0$ such that $G$ is positive semi-definite, then the
minimization problem to obtain $v^t$ in line 3 of
Algorithm~\ref{algo:admm} becomes:
\[
v^t=\arg\min_v \left[g(v) - (\alpha^{t-1}+\rho B^\top A w^{t-1} +
  \rho G
v^{t-1})^\top v + \frac{\rho \beta}{2} v^\top v \right] ,
\]
which may be simpler to solve than the corresponding problem with $G=0$,
since the original quadratic term
$v^\top B^\top B v$ is now replaced by $v^\top v$. 
The additional term $\|w-w^{t-1}\|_H^2$ can play a similar role of preconditioning. 

\begin{algorithm}
\begin{algorithmic}[1]
\STATE Choose $w^0$, $v^0$, and $\alpha^0$
\FOR{ $t =1,2,\ldots$}
\STATE $v^t=\arg\min_v [g(v) - {\alpha^{t-1}}^{\top} B v +
\frac{\rho}{2}\|A {w}^{t-1} - Bv-c\|_2^2 + \frac{\rho}{2} \|v-v^{t-1}\|_G^2]$;
\STATE $w^t=\arg\min_w [\phi(w)  + {\alpha^{t-1}}^{\top} A w +
\frac{\rho}{2}\|Aw-Bv^{t}-c\|_2^2 + \frac{1}{2} \|w-w^{t-1}\|_H^2]$;
\STATE $ \alpha^{t}=\alpha^{t-1} + \rho (Aw^{t}-Bv^t-c)$;
\ENDFOR
\STATE \textbf{Output:} $w^t$, $v^t$, $\alpha^t$.
\end{algorithmic}
\caption{Preconditioned Standard ADMM Algorithm}
\label{algo:admm}
\end{algorithm}

For simplicity, this paper focuses on the scenario that $g(\cdot)$ is strongly convex, and
$\phi(\cdot)$ is smooth.
The results allow $g(\cdot)$ to include a constraint $v \in \Omega$ 
for a convex set $\Omega$ by setting $g(v)=+\infty$ when $v \notin \Omega$.
The same proof technique can also handle other three cases with one objective function
being smooth and one being strongly convex. 

The standard ADMM algorithm assumes that the optimization problem to obtain 
 $w^t$ is simple. If this optimization is difficult to perform, then we
may also consider the linearized ADMM formulation which replaces $\phi(w)$ by
a quadratic approximation $\phi_H(w)$ defined as
\[
\phi_H(w^{t-1};w)= \phi(w^{t-1}) + \nabla \phi(w^{t-1})^\top
(w-w^{t-1}) + \frac12 (w-w^{t-1})^\top H (w-w^{t-1}) .
\]
The resulting algorithm is described in
Algorithm~\ref{algo:lin-admm}. Both $H$ and $G$ are symmetric positive
semi-definite matrices.
By setting $H=\beta' I - \rho A^\top A$, we can replace
the term $w^\top A^\top A w$ by $w^\top w$ in the optimization of line
4 of Algorithm~\ref{algo:lin-admm}. 

\begin{algorithm}
\begin{algorithmic}[1]
\STATE Choose $w^0$, $v^0$, and $\alpha^0$
\FOR{ $t =1,2,\ldots$}
\STATE $v^t=\arg\min_v [g(v) - {\alpha^{t-1}}^{\top} B v +
\frac{\rho}{2}\|A w^{t-1} - B v -c\|_2^2 + \frac{\rho}{2} \|v-v^{t-1}\|_G^2]$;
\STATE $w^t=\arg\min_w [\phi_H(w^{t-1};w) + {\alpha^{t-1}}^{\top} A w + \frac{\rho}{2}\|Aw-Bv^{t}-c\|_2^2]$;
\STATE $ \alpha^{t}=\alpha^{t-1} + \rho (Aw^{t}-Bv^t-c)$;
\ENDFOR
\STATE \textbf{Output:} $w^t$, $v^t$, $\alpha^t$.
\end{algorithmic}
\caption{Preconditioned Linearized ADMM Algorithm}
\label{algo:lin-admm}
\end{algorithm}

This paper compares the convergence behavior of the ADMM algorithm
versus that of the linearized ADMM algorithm for
solving \eqref{eqn:aecco}. Under the assumption that $A$ is
invertible, $g(\cdot)$ is $\lambda$
strongly convex, and $\phi(\cdot)$ is $1/\gamma$ smooth, it is shown that
the standard ADMM achieves a worst case linear convergence rate of
$1/(1+\Theta(\sqrt{\lambda \gamma}))$ 
(with optimally chosen $\rho$) while the linearized ADMM
achieves a slower worst case linear convergence rate of $1/(1+\Theta(\lambda \gamma))$. 

The paper is organized as follows. Section~\ref{sec:related} reviews
related work. Section~\ref{sec:analysis} provides a theoretical analysis for both
standard and linearized ADMM. 
Section~\ref{sec:numerical} provides a simple numerical example to illustrate
the difference in convergence behavior.
Concluding remarks are given in Section~\ref{sec:conclusion}. 

\section{Related Work on ADMM and Linearized ADMM}
\label{sec:related}
In this section, we review some previous work on the convergence
analysis of ADMM and Linearized ADMM, focusing mainly on linear convergence results.

\subsection{Results for ADMM}

Many authors have studied the linear convergence of ADMM in recent
years. For example, the authors in \cite{hong2012linear} presented a
novel proof for the linear convergence of the ADMM algorithm. 
Moreover, the analysis applies for the more general case in which the
object function can be the summation of more than two separable
functions ($\phi$ and $g$ in our case). 
However, the assumption on each separable function is very complex, and
no explicit rate is obtained. Therefore their results are not directly
comparable to ours.

Another work is \cite{deng2012global}, which presented analysis for
the linear convergence of generalized ADMM under certain conditions. 
More comprehensive results for the general form of constraint $A w - B
v=c$ were obtained later in \cite{davis2014faster} using similar
ideas. 
In that paper, they presented an extension of ADMM algorithm called {Relaxed ADMM}, which leads to linear convergence in the following four cases
(it also requires either $A$ or $B$ are invertible):
$\phi$ is strongly convex and smooth;
$g$ is strongly convex, and smooth;
$\phi$ is smooth, and $g$ is strongly convex;
$g$ is smooth, and $\phi$ is strongly convex.
However, their analysis employs a technique for analyzing the dual objective
of ADMM that may be regarded as a Relaxed Peacheman-Rachford
splitting method. It can be used to prove the dual convergence.
In contrast, our analysis uses a very different argument that
can directly bound the convergence of primal objective function and
the duality gap.
Moreover, even when the
required regularity conditions for linear convergence are not
satisfied, our analysis immediately implies a sublinear $1/t$
convergence of duality gap (assuming a finite solution exists
for the underlying problem). Therefore the analysis of this paper
contains a unified treatment that can simultaneously handle both linear and sublinear
convergence depending on the regularity condition. 
In contrast, although sublinear results 
can be obtained using techniques similar to those of \cite{davis2014faster}
(see results in \cite{DavisYin14-sublinear}), they require specialized
treatment and the obtained results are in different forms that are not
compatible with the duality gap convergence of this paper. In this
setting, the operator splitting proof techniques of
\cite{davis2014faster,DavisYin14-sublinear}
and the objective function proof technique of this paper are complementary
to each other. 
Another advantage of our proof technique is that it can be directly
applied to linearized ADMM with minimal modifications.

Our analysis employs a technique similar to  that of
\cite{ouyang2013stochastic} (note that neither linear convergence nor
duality gap convergence was
studied in \cite{ouyang2013stochastic}). 
At the conceptual level, the technique is also closely related to the analysis of
\cite{chambolle2011first}, but the actual execution differs quite
significantly. One may view the analysis of this paper as a refined 
version of those in \cite{ouyang2013stochastic}, in that we simultaneously handle linear
and sublinear cases depending on regularity conditions. 
Moreover, our analysis unifies the techniques used in \cite{ouyang2013stochastic} (which deals with
primal objective convergence)  and the techniques used in \cite{chambolle2011first} (which deals with a special primal-dual objective
convergence); our proof shows that the seemingly different results in these two papers can
be proved using the same underlying argument. Although results similar to ours were
presented in \cite{chambolle2011first} for a procedure related to a
specific form of preconditioned ADMM
(see \cite{chambolle2011first} for discussions), they did not analyze the
standard ADMM (or its linearized version) under the general condition
$Aw - B v =c$.
Therefore results obtained in this paper for ADMM are different from those
of \cite{chambolle2011first}.

Another result on the linear convergence of the standard ADMM can be found
in a recent paper \cite{nishihara2015general}, which uses a different technique than
what's presented in this paper and that of
\cite{deng2012global,davis2014faster}. Their results are not directly comparable
to ours. 
Moreover, some other work on the convergence of ADMM like procedures include
\cite{iutzeler2013explicit,goldfarb2013fast,suzuki2013dual}, 
which focused on different applications that are not related to our
work. 

\subsection{Results for Linearized ADMM}

One advantage of our proof technique is that it also handles
linearized ADMM, with new results not available in
the previous literature. Most of previous work on linearized ADMM does not
consider linear convergence; a few that do consider impose strong
assumptions on the matrices $A, B$, or the functions $f, g$.

There are several papers that considered linear convergence of
Linearized ADMM. For example \cite{hong2012linear} considered
linearized ADMM, but as mentioned earlier, their rate is not explicit
and they impose complex conditions that are incompatible with our
results.
Similarly, a linear convergence result for linearized ADMM was also
obtained in \cite{ling2014decentralized}, but only under the
assumption of $g=0$ and some strong constraints on the matrices $A$
and $B$. Again their results are incompatible with ours.  

Some other work considered Linearized ADMM in the general cases but
without linear convergence. For example, in \cite{ouyang2015accelerated}, the
authors consider the convergence of Linearized ADMM on several
different cases, and obtained sublinear convergence of $1/t$.  
Similar sublinear results can be found in \cite{ouyang2013stochastic} for
stochastic ADMM. As we have pointed out, our proof technique is
closely related to that of \cite{ouyang2013stochastic}, which can
handle both linearized and standard ADMM under the same theoretical framework.

\section{Main Results}
\label{sec:analysis}

This section provides our main results for the standard ADMM and the
 linearized ADMM. 
We will derive upper bounds on their convergence rates, as well as
the worst case matching lower bounds for some specific problems. 

\subsection{Notations}
Given any convex function $h$, we may define its convex conjugate
\[
h^*(\beta) = \sup_{u} [ \beta^\top u - h(u)] ,
\]
and define the Bregman divergence of a convex function $h(u)$ as:
\[
 D_h(u',u)= h(u) -h(u') - \nabla h(u')^\top (u-u') .
\]

We will assume that $\phi$ is $1/\gamma$ smooth: 
\[
\forall w, w', \qquad 
D_\phi(w',w) \leq \frac{1}{2\gamma} \|w'-w\|_2^2 ,
\]
which also implies that
\[
D_\phi(w',w) \geq \frac{\gamma}{2} \|\nabla \phi(w')-\nabla \phi(w)\|_2^2 ,
\qquad
D_{\phi^*}(u',u) \geq \frac{\gamma}{2} \|u'-u\|_2^2 .
\]

We also assume that $g$ is $\lambda$ strongly convex:
\[
\forall v, v', \qquad D_g(v',v) \geq \frac{\lambda}{2} \|v'-v\|_2^2 .
\]

Assume also that $(w_*,v_*,\alpha_*)$ is an optimal solution of
\eqref{eqn:aecco}, which satisfies the equality:
\begin{equation}
A w_* - B v_* -c =0 , 
\quad A^\top \alpha_*= -\nabla \phi(w_*) , 
\quad w_* =\nabla \phi^*(-A^\top \alpha_*), 
\quad B^\top \alpha_* =  \nabla g(v_*) . \label{eq:opt}
\end{equation}

Given any $\alpha$, taking inf over $(w,v)$ with respect to the Lagrangian
\[
\phi(w)+g(v) + \alpha^\top (Aw-Bv-c), 
\]
we obtain the dual
\[
D(\alpha)= - \phi^*(-A^\top \alpha) - g^*(B^\top \alpha) - \alpha^\top
c .
\]
It is clear by definition that for any pair $(w,v)$ that are feasible (that is $Aw -
Bv-c=0$), and any $\alpha$, we have $\phi(w)+g(v) \geq D(\alpha)$. The value
$\phi(w)+g(v) -D(\alpha)$ is referred to as the duality gap. 
Duality gap is always larger than primal suboptimality $[\phi(w)+g(v)]-[\phi(w_*)+g(v_*)]$.
Therefore if the duality gap is zero, then $(w,v)$ solves \eqref{eqn:aecco}. 

We may also introduce the concept of restricted duality gap
as in \cite{chambolle2011first}. 
Consider regions $B_1 \subset \Real^p$, and $B_2 \subset \Real^m$.
Given any $\hat{\alpha}$, $\hat{v}$, we can define the restricted duality gap
\[
\cG_{B_1 \times B_2}(\hat{\alpha},\hat{v})=\sup_{\alpha \in B_1; v \in B_2} \left[ \phi^*(-A^\top \hat{\alpha}) + g(\hat{v}) - \phi^*(-A^\top \alpha) - g(v) 
+ \hat{\alpha}^{\top}(Bv+c) - \alpha^\top (B\hat{v}+c) \right] .
\]
If we pick $(\alpha,v)=(\alpha_*,v_*)$, then
\[
D_{\phi^*}(-A^\top \alpha_*,-A^\top \hat{\alpha}) + D_g(v_*,\hat{v}) =
\phi^*(-A^\top \hat{\alpha}) + g(\hat{v}) - \phi^*(-A^\top \alpha_*) - g(v_*) 
+ \hat{\alpha}^{\top}(Bv_*+c) - \alpha_*^\top (B\hat{v}+c) .
\]
Therefore as long as $(\alpha_*,v_*) \in B_1 \times B_2$, we have
\[
D_{\phi^*}(-A^\top \alpha_*,-A^\top \hat{\alpha}) + D_g(v_*,\hat{v}) 
\leq \cG_{B_1 \times B_2}(\hat{\alpha},\hat{v})  .
\]
Assume $AA^\top$ is invertible, and let 
\begin{equation}
A^{+}=A^\top (A A^\top)^{-1} \label{eq:Ainv}
\end{equation}
be the pseudo-inverse of $A$, then we may let $\hat{w}=A^{+}(B\hat{v}+c)$. It follows that $A\hat{w}-B\hat{v}-c=0$.
If we set $B_1 \times B_2 = \Real^p \times \Real^m$, then we recover the
unrestricted duality gap:
\[
\cG_{\Real^p \times \Real^m}(\hat{\alpha},\hat{v})  = [\phi(\hat{w}) + g(\hat{v})]
- D(\hat{\alpha}) ,
\]
where the maximum over $(\alpha,v)$ is taken at $-A^\top\alpha=\nabla
\phi(\hat{w})$ and $v=\nabla g^*(B^\top\hat{\alpha})$.

\subsection{Standard ADMM}

In general, we have the following result.
\begin{theorem}
  Assume that $\phi$ is $1/\gamma$ smooth and $g$ is $\lambda$
  strongly convex.   Assume that we can write $H=A^\top \tilde{H} A$.
  Let $\sMH$ and $\stH$ be the largest eigenvalues of $H$ and $\tilde{H}$
  respectively,
 $\sA$ be be the smallest eigenvalue value of $(AA^\top)^{1/2}$,
  $\sB$ be the largest singular value of $B$, $\sG$ be the largest
  singular value of $G$. Consider $s \in [0,1)$ and $\theta>0$ such that
  \[
  \theta \leq \min \left(\frac{\gamma \rho \sA^2}{\gamma\sMH+1} ,\
    \frac{s \rho}{\stH},
  \frac{(1-s)\lambda}{(\rho+\stH) \sB^2+(1-s)\rho \sG}\right) .
  \]
  Let
  $\tilde{\alpha}^t=\alpha^t+ \tilde{H}A(w^t-w^{t-1})$.
  Then for all $(\alpha,v)$ and $w=\nabla \phi^*(-A^\top \alpha)$, Algorithm~\ref{algo:admm} produces approximate solutions that
  satisfy 
  \begin{align}
\sum_{t=1}^T
(1+\theta)^{t-T} r_t 
\leq & (1+\theta)^{-T} \delta_0 - \delta_T ,  \label{eq:admm-primal}\\
\sum_{t=1}^T (1+\theta)^{t-T} r_t^*
\leq& (1+\theta)^{-T} \delta_0 - \delta_T , \label{eq:admm-dual}
\end{align}
where
\begin{align*}
r_t =& \phi(w^{t}) + g(v^t) - \phi(w) - g(v)
- \tilde{\alpha}^{t\top}(Aw-Bv-c)
+ \alpha^\top (Aw^{t}-Bv^t-c) , \\
r_t^*=&  \phi^*(-A^\top \tilde{\alpha}^{t}) + g(v^t) - \phi^*(-A^\top \alpha) - g(v) 
+ \tilde{\alpha}^{t\top}(Bv+c)- \alpha^\top (Bv^t+c) , \\
\delta_t=& \frac{\rho}{2} \|Aw^{t}-Bv-c\|_2^2 +\frac12 \|Aw^{t}-Bv-c\|_{\tilde{H}}^2
+ \frac{\rho (1+\theta)}{2} \|v^t-v\|_G^2   
+ \frac{1+\theta}{2\rho}  \|\alpha-\alpha^{t}\|_2^2 .
\end{align*}
\label{thm:admm}
\end{theorem}
For arbitrary $(\alpha,v)$, the left hand side of \eqref{eq:admm-primal} and \eqref{eq:admm-dual} 
can be difficult to understand. We may choose specific values of $(\alpha,v)$ so that
the results are easier to interpret. 
By setting $(\alpha,w,v)=(\alpha_*,w_*,v_*)$ in
Theorem~\ref{thm:admm}, and using \eqref{eq:opt}, we obtain the
following corollary.
\begin{corollary}
Under the conditions of Theorem~\ref{thm:admm}, we have
  \begin{align*}
&\sum_{t=1}^T
(1+\theta)^{t-T} \left[\max(D_\phi(w_*,w^{t}),D_{\phi^*}(-A^\top \alpha_*,-A^\top \tilde{\alpha}^{t})) + D_g(v_*,v^t) \right]\\
&+ \frac{\rho}{2} \|A(w^T-w_*)\|_2^2 
+ \frac12 \|A(w^{T}-w_*)\|_{\tilde{H}}^2 
+ \frac{1+\theta}{2\rho} \|\alpha^T-\alpha_*\|_2^2
+ \frac{\rho(1+\theta)}{2} \|v^T-v_*\|_G^2\\
\leq &  
\frac{(1+\theta)^{-T}}{2}
\left[
 {\rho} \|A(w^{0}-w_*)\|_2^2 
+ \|A(w^{0}-w_*)\|_{\tilde{H}}^2
+ \frac{1+\theta}{\rho}  \|\alpha^{0}-\alpha_*\|_2^2 
+ \rho (1+\theta) \|v^{0}-v_*\|_G^2 
\right] .
\end{align*}
\label{cor:admm-bregman}
\end{corollary}
Using the definition of restricted duality gap, 
it is easy to see that \eqref{eq:admm-dual} directly implies 
an upper bound of restricted duality gap, which is the same style as results
of \cite{chambolle2011first}. 
Our result is more general than
those of \cite{chambolle2011first} because the results can also be
expressed in the form of Corollary~\ref{cor:admm-bregman}, as well as
in terms of unrestricted duality gap, as stated below.
\begin{corollary}
  Under the conditions of Theorem~\ref{thm:admm}, and let $A^{+}$ be
  the psudo-inverse of $A$. Define
  \begin{align*}
    \delta_*^0=& 
    \left[
      (\rho+\stH) \|A(w^{0}-w_*)\|_2^2 
      + \frac{1+\theta}{\rho}  \|\alpha^{0}-\alpha_*\|_2^2 
      + \rho (1+\theta) \|v^{0}-v_*\|_G^2 
    \right] \\
  b_1(\delta) =& \sup_{u}
  \left\{
    \|\alpha^0+ A^{-T}\nabla \phi(A^{+}(Bu+c)\|_2^2
      : \|u-v_*\|_2^2 \leq \frac{\delta}{ \rho(1+\theta)}
    \right\} \\
  b_2(\delta) =& \sup_{\beta}
  \left\{
    \|v^0-\nabla g^*(B^\top \beta)\|_2^2: 
      \|\beta-\beta_*\|_2 \leq   \sqrt{ \frac {\rho\delta}{1+\theta}}
      + \stH \sqrt{\frac{2(2+\theta)\delta}{\rho}}
    \right\} \\
    b(\delta)=&\frac{1+\theta}{2\rho} b_1(\delta)+ \frac{1+\theta}{2/\rho} \sG b_2(\delta) 
    +  (\rho+\stH) (\sB^2 b_2(\delta) 
    + \|Aw^{0}-Bv^0-c\|_2^2) .
\end{align*}
Then we have the following bound in duality gap
\[
[\phi(A^{+}(Bv^T+c)) +  g(v^T)] - D(\tilde{\alpha}^T) 
\leq
(1+\theta)^{-T}   b((1+\theta)^{-T} \delta_*^0) .
\]
Moreover, define
\[
\bar{v}^T = \frac{\sum_{t=1}^T (1+\theta)^{t} v^t}{\sum_{t=1}^T
  (1+\theta)^{t}} ,
\qquad
\bar{\alpha}^T = \frac{\sum_{t=1}^T (1+\theta)^{t}
  \tilde{\alpha}^t}{\sum_{t=1}^T (1+\theta)^{t}} .
\]
Then
\[
[\phi(A^{+}(B\bar{v}^T+c)) +  g(\bar{v}^T)] - D(\bar{\alpha}^T) 
\leq
\frac{
b(\delta_*^0)}
{\sum_{t=1}^T (1+\theta)^{t}} .
\]
\label{cor:admm-duality-gap}
\end{corollary}

In the above results,  we consider the simple case of $H=0$. Then the optimal value of $\theta$ is achieved when we take
\[
\rho = \frac{\sqrt{\sB^2+\sG}}{\sA} \sqrt{\frac{\lambda}{\gamma}} , \qquad \theta=
\sA \sqrt{(\sB^2+\sG)\gamma \lambda} .
\]
When $\theta>0$, this implies the following convergence  from Corollary~\ref{cor:admm-bregman}:
\begin{align*}
&\max[D_\phi(w_*,w^T), D_{\phi^*}(-A^\top \alpha_*,-A^\top \tilde{\alpha}^T)]
+  D_g(v_*,v^T) \\
& \quad + \frac{\rho}{2} \|A(w^T-w_*)\|_2^2
+ \frac{1}{2\rho} \|\alpha^T-\alpha_*\|_2^2
+ \frac{\rho}{2} \|v^T-v_*\|_G^2\\
\leq &
\left(1+ \sA \sqrt{(\sB^2+\sG)\gamma \lambda}\right)^{1-T}
\left[
  \frac{\rho}{2} \|A(w^{0}-w_*)\|_2^2 
+ \frac{1}{2\rho}  \|\alpha^{0}-\alpha_*\|_2^2 
+ \frac{\rho}{2}  \|v^{0}-v_*\|_G^2 
\right] .
\end{align*}
This implies $\|w_*-w^T\|_2 = O((1+\theta)^{-T})$,
$\|v_*-v^T\|_2=O((1+\theta)^{-T})$, and $\|\alpha_*-\alpha^T\|_2 = O((1+\theta)^{-T})$.

The linear convergence result holds when $\theta>0$. However, even
when $\theta=0$ (and $H \neq 0$), we can still obtain the following sublinear
convergence from Corollary~\ref{cor:admm-bregman}:
\begin{align*}
& \max\left[D_\phi(w_*,\bar{w}^T), \frac{1}{T}\sum_{t=1}^T D_{\phi^*}(-A^\top \alpha_*,-A^\top \tilde{\alpha}^t)\right]
+  D_g(v_*,\bar{v}^T) \\
\leq &
\frac{1}{2T}
\left[
  (\rho+\stH) \|A(w^{0}-w_*)\|_2^2 
  + \frac{1}{\rho} \|\alpha^{0}-\alpha_*\|_2^2
  + \rho \|v^{0}-v_*\|_2^2
\right] ,
\end{align*}
where
$\bar{w}^T= T^{-1}\sum_{t=1}^T w^t$,
$\bar{v}^T= T^{-1}\sum_{t=1}^T v^t$.
This result does not require any assumption on $\phi$, $g$, $A$,
$B$. 

Similar results hold for unrestricted duality gap under the conditions of
Corollary~\ref{cor:admm-duality-gap}. For example, when $\theta=0$,
but $AA^\top$ is invertible, we
obtain the sublinear convergence of duality-gap below. 
\[
[\phi(A^{+}(B\bar{v}^T+c)) +  g(\bar{v}^T)] - D(\bar{\alpha}^T) 
\leq \frac{b(\delta_*^0)}{T} .
\]
This bound can be compared to
the main result of \cite{chambolle2011first} stated in terms of the 
restricted duality gap
(in which the authors studied a method that is related to, but not
identical to ADMM). Their result did not imply  a bound on the unrestricted
duality gap because they did not obtain a counterpart of Corollary~\ref{cor:admm-bregman}.

In the case of $\phi$ being smooth but $g$ is not a strongly convex function, given any
$\epsilon>0$, we can set $\lambda=\epsilon$, and apply ADMM with
$g(v)$ replaced by the strongly convex function $g(v)+\lambda v^\top v$.
With $\rho$ chosen optimally, this leads to
\[
\phi(A^{+}(Bv^T+c)) + g(v^T) - D(\alpha^T) = O(\epsilon)
\]
when we take $T= \ln (1/(\gamma\epsilon))/\sqrt{\gamma\epsilon}$. 

\subsection{Linearized ADMM}

For Linearized ADMM, we have the following counterpart of
Theorem~\ref{thm:admm}.
Here we need to assume that $A$ is invertible and $H$ is sufficiently
large so that $\smH \geq \gamma^{-1}$. 
\begin{theorem}
  Assume that $\phi$ is $1/\gamma$ smooth and $g$ is $\lambda$
  strongly convex, and $A$ is a square invertible matrix.  
  Assume that we can write $H=A^\top \tilde{H} A$.
  Let $\smH$ and $\stH$ be the smallest eigenvalue of $H$ and 
  the largest eigenvalue of $\tilde{H}$  respectively, and we
  assume that $\smH \geq \gamma^{-1}$. 
  Let $\sA$ be the smallest eigenvalue value of $(AA^\top)^{1/2}$,
  $\sB$ be the largest singular value of $B$, $\sG$ be the largest
  singular value of $G$. Consider $s \in [0,1)$ and $\theta>0$ such that
  \[
  \theta \leq \min\left(
    \frac{\rho \sA^2}{\smH},
       \frac{s \rho}{\stH},
  \frac{(1-s)\lambda}{(\rho+\stH) \sB^2+(1-s)\rho\sG}\right) .
  \]
  Let $\tilde{\alpha}^t=
  \alpha^t + \tilde{H}A( w^t-w^{t-1})$.
  Then for any $(\alpha,v)$, and $w =\nabla \phi^*(-A^\top \alpha)$,
  Algorithm~\ref{algo:lin-admm} produces approximate solutions that
  satisfy
  \begin{align}
\sum_{t=1}^T (1+\theta)^{t-T} r_t 
\leq& (1+\theta)^{-T} \delta_0 - \delta_T ,
  \label{eq:lin-admm-primal} \\
\sum_{t=1}^T (1+\theta)^{t-T} r_t^*
\leq& (1+\theta)^{-T} \delta_0 - \delta_T ,
\label{eq:lin-admm-dual}
\end{align}
where
\begin{align*}
{r}_t =&
\phi(w^{t-1}) + g(v^t) - \phi(w) - g(v)
- \tilde{\alpha}^{t\top}(Aw- Bv-c) + \alpha^\top
(Aw^{t-1}-Bv^t-c) ,\\
{r}_t^* =& 
\phi^*(-A^\top \tilde{\alpha}^{t}) + g(v^t) - \phi^*(-A^\top \alpha) - g(v)
+ \tilde{\alpha}^{t\top} (Bv+c) - \alpha^\top (Bv^t+c) , \\
\delta_t=&
\frac12
\left[  \|Aw^t-Bv-c\|_{\tilde{H}}^2 
+ \rho \|Aw^t-Bv-c\|_2^2
+ \rho (1+\theta)\|v^t-v\|_G^2
+ \frac{1+\theta}{\rho} \|\alpha^t-\alpha\|_2^2
\right] .
\end{align*}
  \label{thm:lin-admm}
\end{theorem}
Similar to the corollaries of Theorem~\ref{thm:admm}, we have the
following three corollaries of Theorem~\ref{thm:lin-admm}.
\begin{corollary}
Under the conditions of Theorem~\ref{thm:lin-admm}, we have
 \begin{align*}
\sum_{t=1}^T
& (1+\theta)^{t-T} \left[\max(D_\phi(w_*,w^{t-1}),D_{\phi^*}(-A^\top
  \alpha_*,-A^\top \tilde{\alpha}^{t})) + D_g(v_*,v^t) \right]\\
&+ \frac{1}{2} \|w^T-w_*\|_H^2 
+ \frac{\rho}{2} \|A(w^T-w_*)\|_2^2 
+ \frac{1+\theta}{2\rho} \|\alpha^T-\alpha_*\|_2^2
+ \frac{\rho(1+\theta)}{2} \|v^T-v_*\|_G^2\\
\leq& 
\frac{(1+\theta)^{-T}}{2}
\left[
 (\rho + \stH) \|A(w^{0}-w_*)\|_2^2 
+ \frac{1+\theta}{\rho}  \|\alpha^{0}-\alpha_*\|_2^2 
+ \rho (1+\theta) \|v^{0}-v_*\|_G^2 
\right] .
\end{align*}
\end{corollary}

\begin{corollary}
  Under the conditions of Theorem~\ref{thm:admm}, and let $A^{+}$ be
  the psudo-inverse of $A$.  If we define
  $\delta_*^0$, $b(\delta)$, $\bar{v}^T$ and $\bar{\alpha}^T$
  as in Corollary~\ref{cor:admm-duality-gap}, then
  \begin{align*}
  [\phi(A^{+}(Bv^T+c)) +  g(v^T)] - D(\tilde{\alpha}^T) 
  \leq&
  (1+\theta)^{-T}   b((1+\theta)^{-T} \delta_*^0) , \\
  [\phi(A^{+}(B\bar{v}^T+c)) +  g(\bar{v}^T)] - D(\bar{\alpha}^T) 
  \leq&
  \frac{
    b(\delta_*^0)}
  {\sum_{t=1}^T (1+\theta)^{t}} .
\end{align*}
\end{corollary}

The requirement of $\smH \geq \gamma^{-1}$ is the key difference
between Theorem~\ref{thm:admm} and Theorem~\ref{thm:lin-admm}.
The fast convergence of ADMM requires that $H$ to be of order
$O(\rho)$, which may be smaller than $\Theta(\gamma^{-1})$.
Consider the case that $H = \Theta (\gamma^{-1} I)$ for linearized
ADMM, then the optimal $\rho$ can be chosen as $\rho = \Theta (\gamma^{-1})$. 
This leads to a linear convergence with $\theta = \Theta(\lambda \gamma)$.
The rate is slower than that of the standard ADMM, which can achieve $\theta=\Theta(\sqrt{\lambda
  \gamma})$ at the optimal choice of $\rho$. 

Similar to the case of standard ADMM, we could take $\theta=0$: as
long $\phi$ is $1/\gamma$ smooth, and $H$ satisfies $\smH \geq
2/\gamma$, we can achieve the following sublinear convergence without
additional assumptions: 
 \begin{align*}
& \sum_{t=1}^T
\left[\max(D_\phi(w_*,w^{t-1}),D_{\phi^*}(-A^\top\alpha_*,-A^\top \tilde{\alpha}^t)) + D_g(v_*,v^{t})\right]\\
\leq& 
\frac{1}{2}
\left[
(\rho+\stH) \|A (w^{0}-w_*)\|_2^2 
+ \rho \|v^0-v_*\|_G^2 
+ \rho^{-1} \|\alpha^0-\alpha_*\|_2^2
\right] .
\end{align*}
A similar result holds for duality gap convergence when $A$ is a square invertible
matrix. 

\subsection{Lower Bounds}

We consider the quadratic case that $A=B=I$, $c=0$, and 
\[
\phi(w) = \frac{1}{2} w^\top Q w , \qquad g(v) = \frac{1}{2}
v^\top \Lambda v .
\]
The optimal solution is
\[
w_* = v_*= \alpha_* = 0 .
\]
We show that with appropriately chosen $Q$ and $\Lambda$ so that 
$Q$ is $1/\gamma$ smooth, and both $\Lambda$ and $Q$ are $\lambda$
strongly convex, the convergence rate of ADMM can be 
$1-\Theta(\sqrt{\gamma\lambda})$ and the convergence 
rate of linearized ADMM can be $1-\Theta(\gamma \lambda)$. 

\subsubsection*{ADMM}
We assume that $Q$ and $\Lambda$ are diagonal matrices.

The ADMM iterate satisfies the following equations (with $G=0$):
\begin{align*}
v^t =& (\Lambda+\rho I)^{-1} (\alpha^{t-1} + \rho w^{t-1}) \\
w^t =& (Q + \rho I)^{-1} ( \rho v^t - \alpha^{t-1} ) \\
\alpha^t=& \alpha^{t-1} + \rho(w^t - v^t) ,
\end{align*}
which implies
\begin{align*}
v^t =& (\Lambda+\rho I)^{-1} (\alpha^{t-1} + \rho w^{t-1}) \\
w^t =& (\Lambda+\rho I)^{-1} (Q + \rho I)^{-1} ( \rho^2 w^{t-1} - \Lambda \alpha^{t-1} ) \\
\alpha^t=& (\Lambda+\rho I)^{-1}(Q+\rho I)^{-1} Q (\Lambda 
\alpha^{t-1}-\rho^2 w^{t-1} ) .
\end{align*}

We may write $[w^t;\alpha^t]= M [w^{t-1};\alpha^{t-1}]$. 
Now we take $Q=\Lambda=\mathrm{diag}(\lambda,1/\gamma)$,
where we assume that $\lambda \leq 1/\gamma$. 
Then the largest eigenvalue of $M$, which determines the rate of
convergence of ADMM, is
\[
\max\left[\frac{\rho^2 + \lambda^2}{(\rho+\lambda)^2} ,
\frac{\rho^2 \gamma^2 + 1}{(\rho\gamma+1)^2} \right] .
\]
The optimal $\rho$ to minimize the above is
$\rho=\sqrt{\lambda/\gamma}$, and the maximum value is $(1+\gamma
\lambda)/(1+\sqrt{\gamma \lambda})^2$. 
This special case matches the convergence rate behavior of
$1-\Theta(\sqrt{\gamma \lambda})$ we proved for the ADMM method.

\subsubsection*{Linearized ADMM}
We assume that $H$, $Q$, and $\Lambda$ are diagonal matrices. 
The linearized ADMM iterate satisfies the following equations (with $G=0$):
\begin{align*}
v^t =& (\Lambda+\rho I)^{-1} (\alpha^{t-1} + \rho w^{t-1}) \\
w^t =& (H + \rho I)^{-1} ((H-Q)w^{t-1}+ \rho v^t - \alpha^{t-1} ) \\
\alpha^t=& \alpha^{t-1} + \rho(w^t - v^t) ,
\end{align*}
which implies that 
\begin{align*}
v^t =& (\Lambda+\rho I)^{-1} (\alpha^{t-1} + \rho w^{t-1}) \\
w^t =& (\Lambda+\rho I )^{-1} (H + \rho I)^{-1} (((\rho I+\Lambda)(H-Q) + \rho^2 I)w^{t-1}
- \Lambda \alpha^{t-1} ) \\
\alpha^t=& 
(\Lambda+\rho I)^{-1}(H+\rho I)^{-1} H (\Lambda 
\alpha^{t-1} + \rho (\Lambda - (\rho I+\Lambda) H^{-1}Q) w^{t-1} ) .
\end{align*}
Now let $\lambda \leq 1/\gamma$, and we take $Q=\Lambda=\mathrm{diag}(\lambda,1/\gamma)$, and $H = \mathrm{diag}(2/\gamma,2/\gamma)$.
It follows that the convergence rate of linearized ADMM is no faster
than the largest eigenvalue of
\[
M= \frac{1}{(\rho +h) (\rho +\lambda)} 
\left[
\begin{array}{cc}
\rho^2 + (\rho + \lambda)(h-q)& -\lambda \\
\rho(\lambda h - (\rho+\lambda)q) & \lambda h 
\end{array}
\right]
\]
with $q=\lambda$ and $h=2/\gamma$. When $\rho \leq h-\lambda$, 
the largest eigenvalue of $M$ is no less than
\[
\frac{\rho^2 + (h-q) \rho + (h-q)\lambda}{\rho^2 + (h+\lambda) \rho + h
  \lambda} \geq \frac{h-\lambda}{h+\lambda} = 1 - O(\lambda\gamma).
\]
Similarly, it is also not difficult to check that the eigenvalue is no less
than $1-O(\lambda\gamma)$ when $\rho \geq h-\lambda$. 
It follows that this special case matches the convergence rate behavior of
$1-\Theta(\gamma \lambda)$ we proved for the linearized ADMM method.

\section{Numerical Illustration}
\label{sec:numerical}

Although we have obtained both the worst case upper bounds and matching lower bounds
for ADMM and Linearized ADMM. The analysis shows that in the worst case
ADMM converges at a faster rate of 
$1-\Theta(\sqrt{\lambda \rho})$
while in the worst case Linearized ADMM converges at a slower rate of
$1-\Theta(\lambda \rho)$. 

However, for any specific problem, both methods can converge faster
than the corresponding worst case upper bounds obtained in this
paper. In this section, we use a simple example to
illustrate the real convergence behavior of ADMM versus linearized ADMM
methods at different choices of $\rho$'s, to illustrate the phenomenon 
that the former can converge significantly faster than the latter.

Consider the following 1-dimensional problem:
\[
\phi(w)= \frac w{\sqrt\gamma}\arctan(\frac w{\sqrt\gamma}) -
\frac12 \ln(1+\frac{w^2}\gamma) + \frac{\mu}{2} w^2,
\quad
g(v) =  \frac1{12} v^4+\frac{\lambda}{2} v^2.
\]
with $A=B=I$ and $c=0$. 
It can be checked that $\phi(w)$ is $1/\gamma+\mu$ smooth and $\mu$
strongly convex; $g(v)$ is $\lambda$-strongly convex. 

 We compare the convergence of ADMM versus linearized ADMM with
different values of $\rho$.  In linearized ADMM, and we set $h=2(\mu+1/\gamma)$.
Note that for this problem, $w_*=v_*=0$, and we can define the error of a
solution $(w,v)$ as $\sqrt{w^2+v^2}$. 

\begin{figure}[h]
\centering
\subfigure[$\rho=0.2\sqrt{\lambda/\gamma}$]{\includegraphics[width=0.32\textwidth]{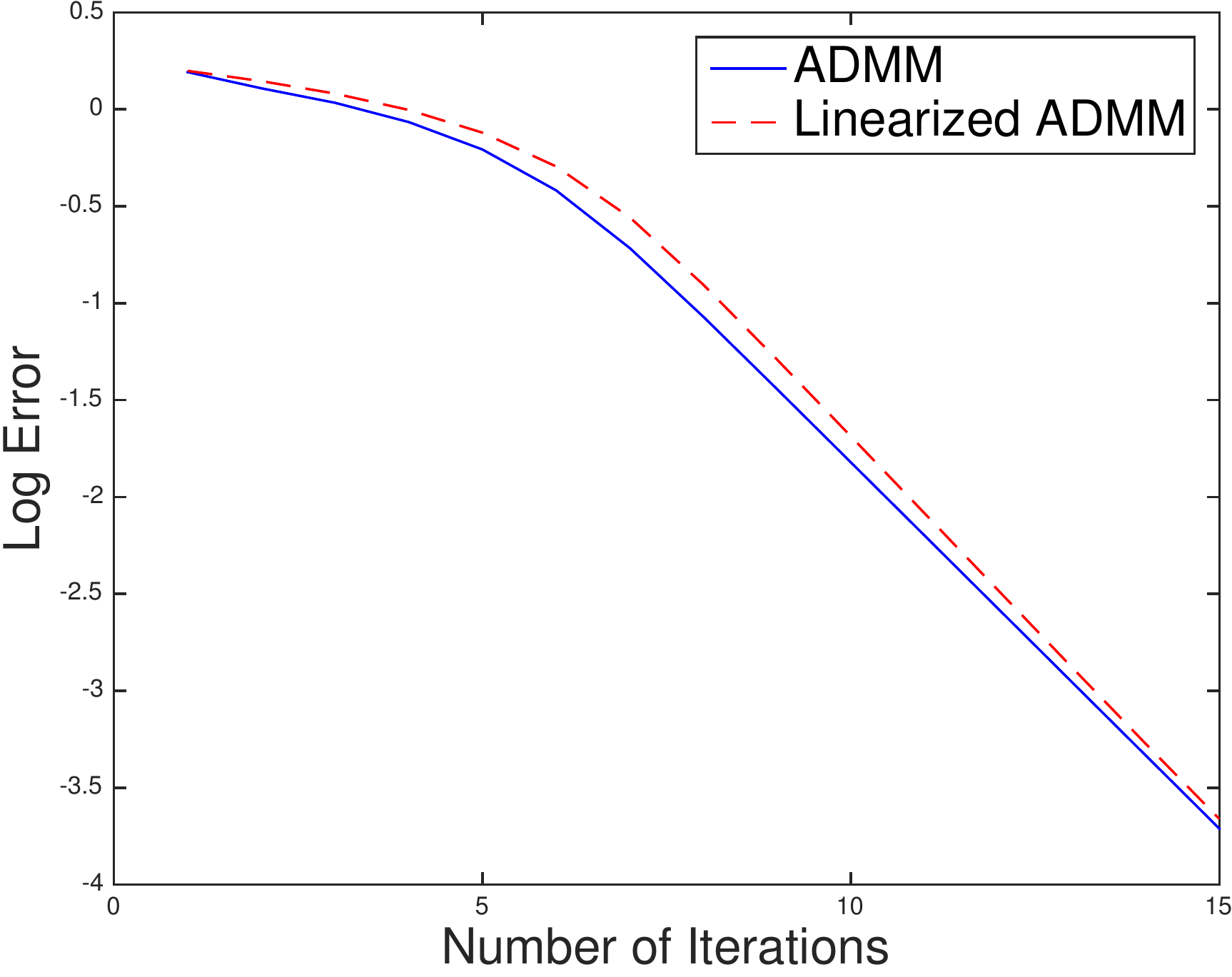}}
\subfigure[$\rho=\sqrt{\lambda/\gamma}$]{\includegraphics[width=0.32\textwidth]{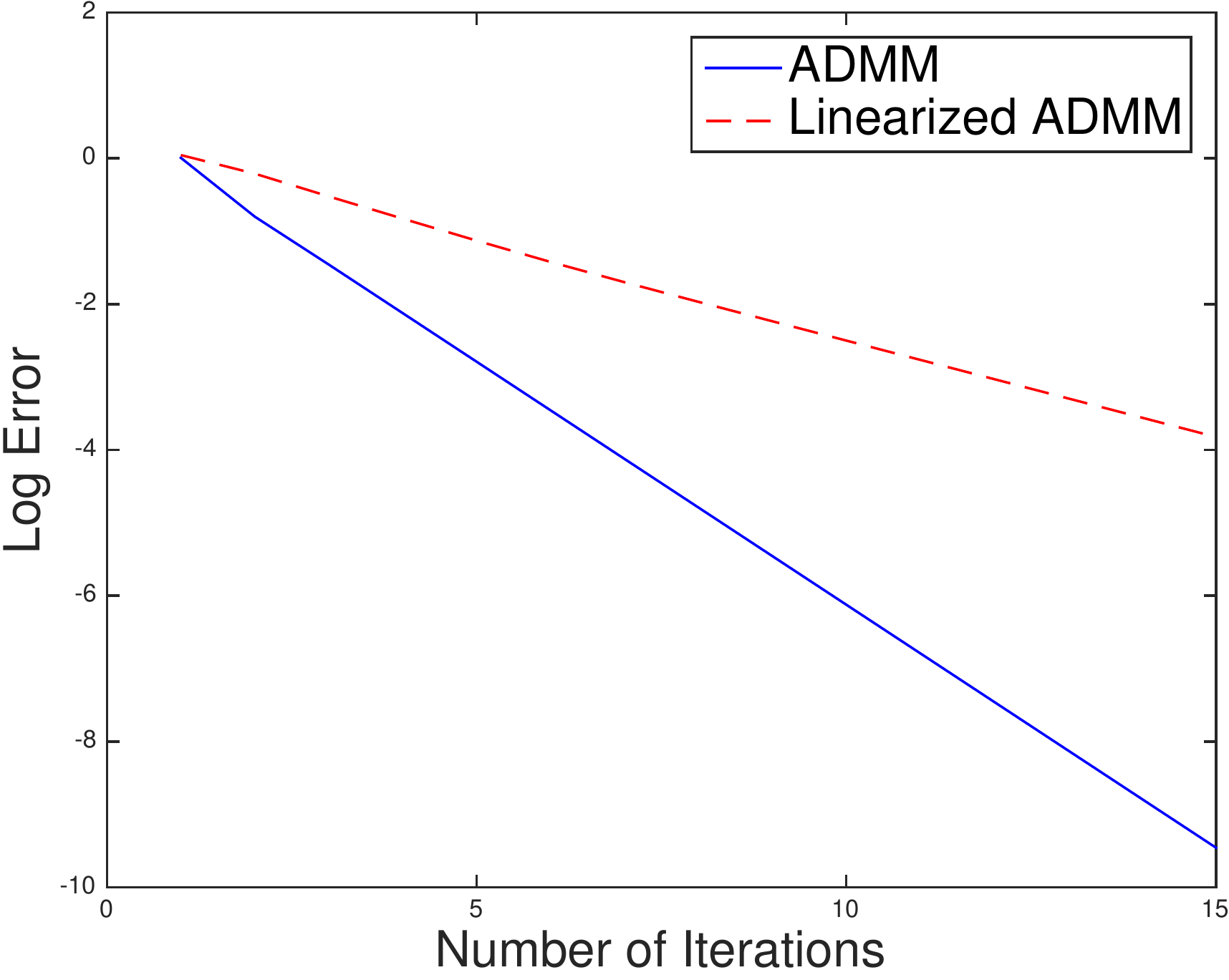}}
\subfigure[$\rho=5\sqrt{\lambda/\gamma}$]{\includegraphics[width=0.32\textwidth]{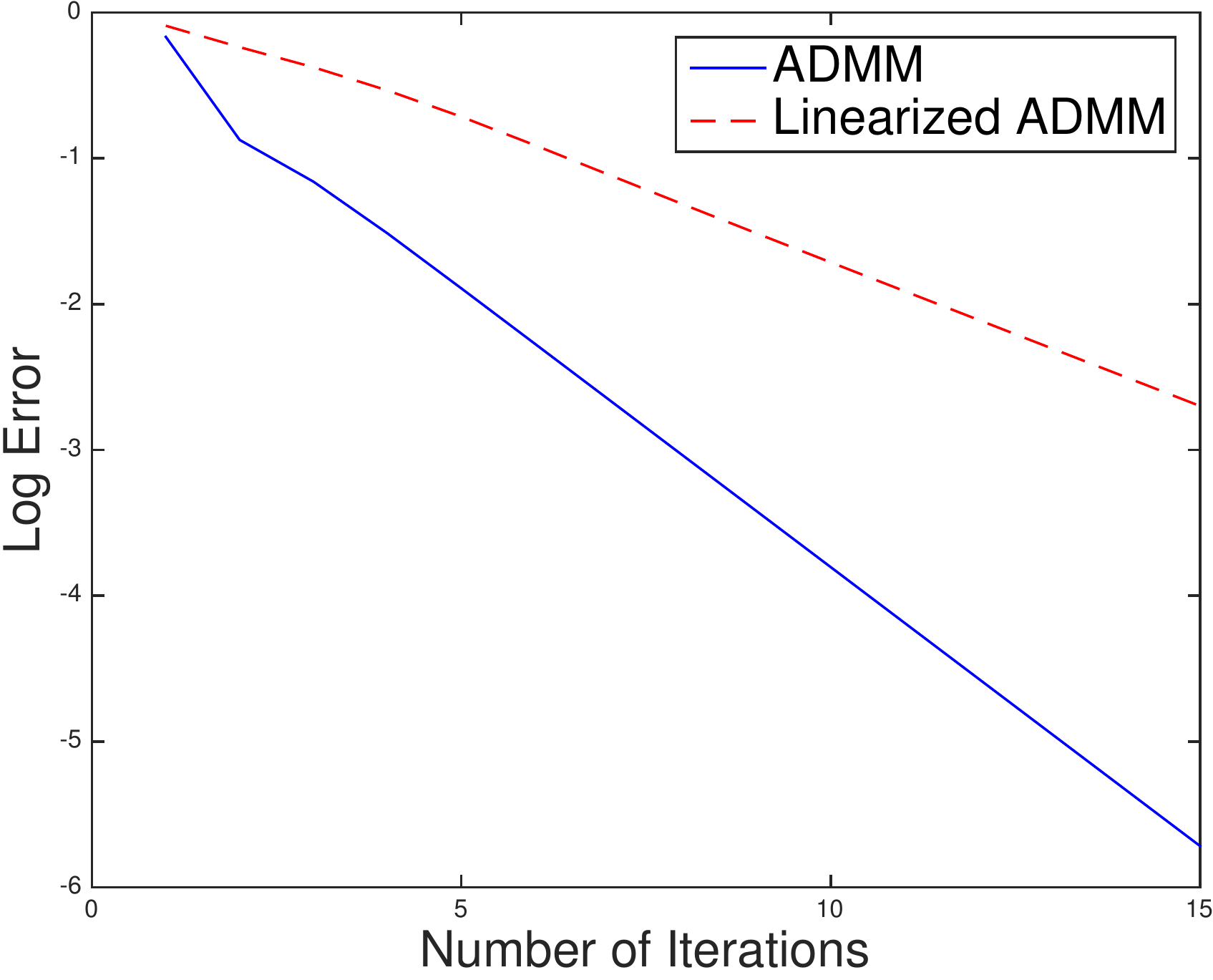}}
\caption{Convergence of ADMM versus that of Linearized ADMM ($\gamma=0.1$,$\lambda=\mu=0.2$)}
\label{fig:conv1}
\end{figure}

Figure~\ref{fig:conv1} shows the convergence behavior when
$\gamma=0.1$, and $\lambda=\mu=0.2$. This is the situation that
$\lambda \gamma=0.02$ is relatively small. In this case,
we compare three different values of
$\rho$'s: $\rho= 0.2 \sqrt{\lambda/\gamma}$,
$\rho=\sqrt{\lambda/\gamma}$, and $\rho=5\sqrt{\lambda/\gamma}$. 
The corresponding convergence rates for ADMM are $0.51$, $0.21$, and $0.41$;
the corresponding convergence rates for linearized ADMM are $0.51$,
$0.53$, and $0.64$.
This shows that ADMM is superior to Linearized ADMM for $\rho$'s. 
Moreover, it achieves relatively fast convergence rate at the optimal
choice of $\rho=\sqrt{\lambda/\gamma}$, while 
Linearized ADMM is relatively insensitive to $\rho$. 

\begin{figure}[h]
\centering
\subfigure[$\rho=0.2\sqrt{\lambda/\gamma}$]{\includegraphics[width=0.32\textwidth]{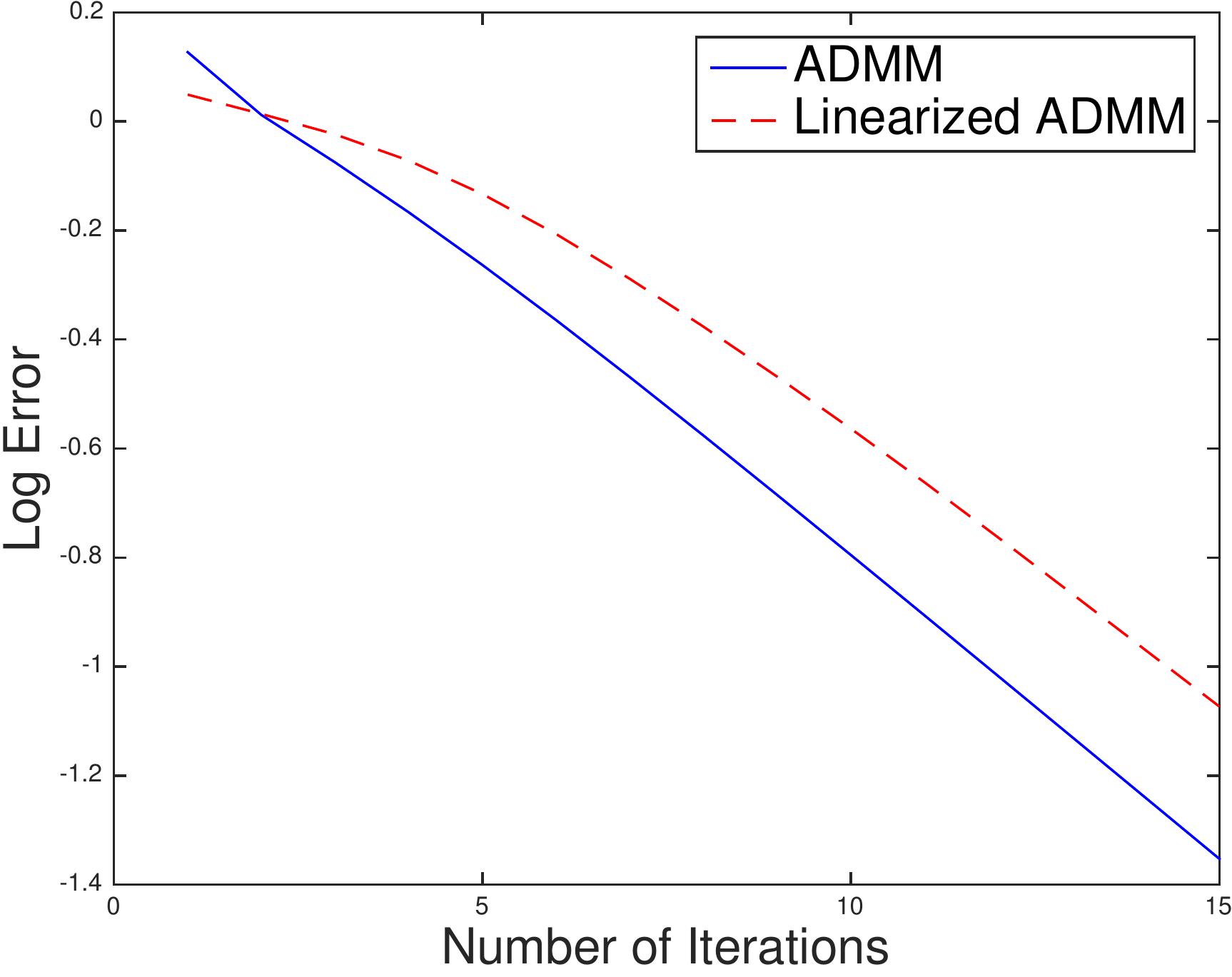}}
\subfigure[$\rho=\sqrt{\lambda/\gamma}$]{\includegraphics[width=0.32\textwidth]{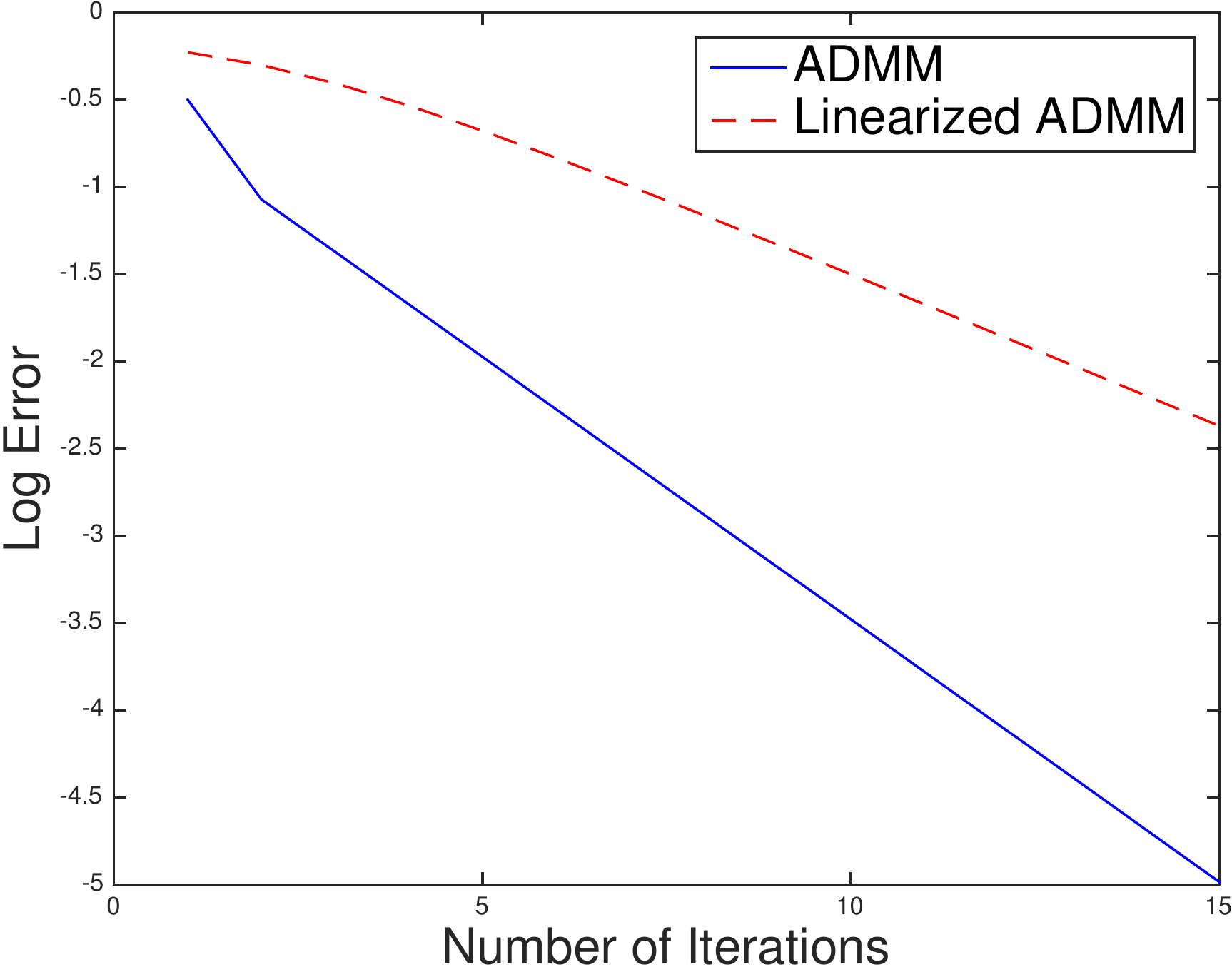}}
\subfigure[$\rho=5\sqrt{\lambda/\gamma}$]{\includegraphics[width=0.32\textwidth]{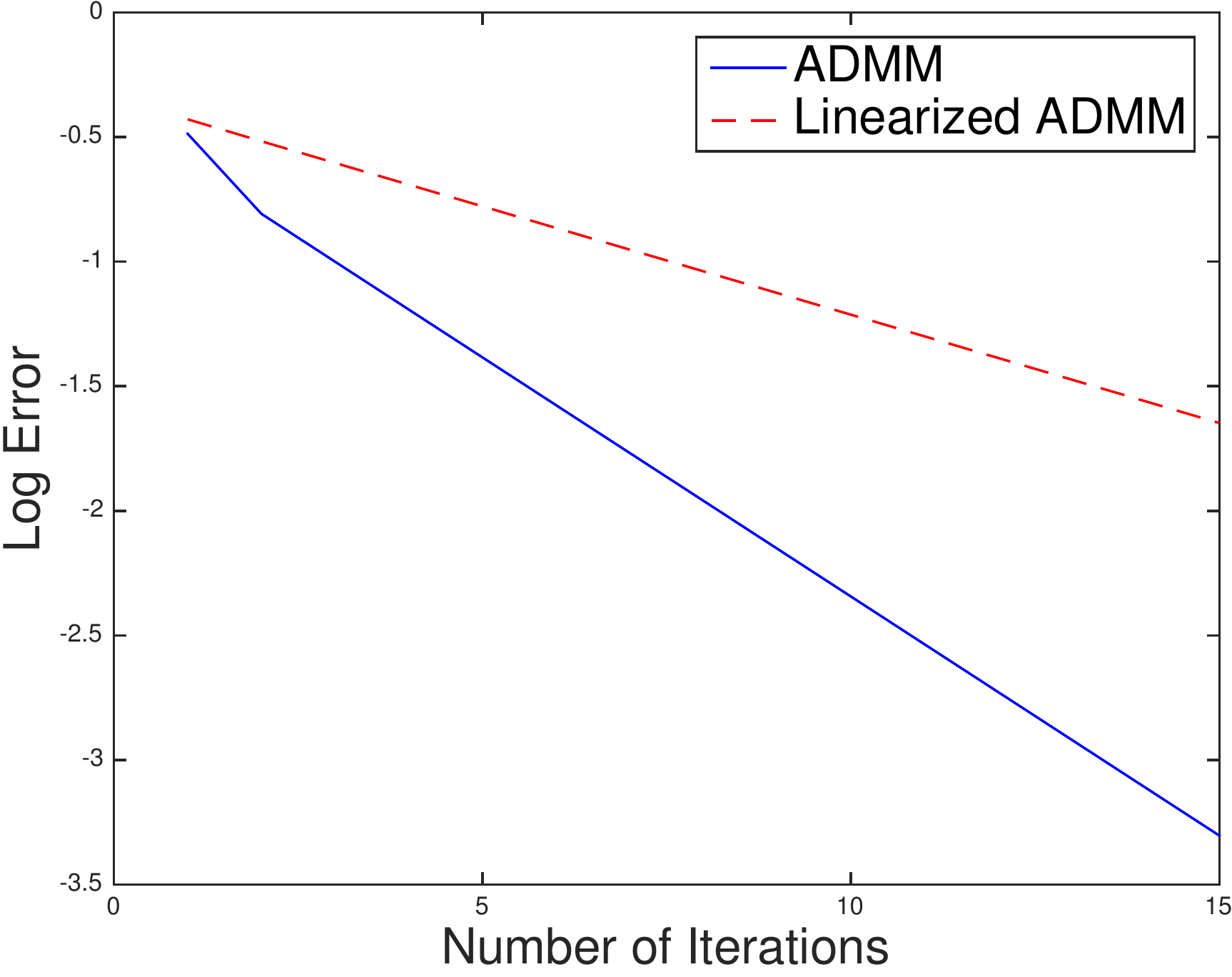}}
\caption{Convergence of ADMM versus that of Linearized ADMM ($\gamma=\lambda=\mu=1$)}
\label{fig:conv2}
\end{figure}

Figure~\ref{fig:conv2} shows the convergence behavior when
$\gamma=\lambda=\mu=1$. This is the situation that
 $\lambda\gamma=1$ is relatively large. 
We compare three different values of
$\rho$'s: $\rho= 0.2 \sqrt{\lambda/\gamma}$,
$\rho=\sqrt{\lambda/\gamma}$, and $\rho=5\sqrt{\lambda/\gamma}$. 
the corresponding convergence rates for ADMM are $0.78$, $0.49$, and $0.64$;
the corresponding convergence rates for linearized ADMM are $0.82$,
$0.69$, and $0.82$.
The relatively convergence behaviors of ADMM and linearized ADMM are consistent
with those of Figure~\ref{fig:conv1}.

\section{Conclusion}
\label{sec:conclusion}

This paper presents a new duality gap convergence analysis of standard ADMM versus linearized
ADMM
under conditions commonly studied in the literature.
It is shown that in the worst case, the standard ADMM converges with an accelerated rate
that is faster than that of the linearized ADMM. 
Matching lower bounds are obtained for specific problems. 
A simple numerical example 
illustrates this behavior. One consequence of our analysis is that the
standard
ADMM does not require Nesterov's acceleration scheme in theory because it
already
enjoys the squared root convergence rate for smooth-strongly convex problems.
On the other hand, linearized ADMM may still benefit from extra acceleration steps.
Finally the results obtained in this paper only show the worst case
behaviors for both algorithms (under appropriate assumptions commonly
used in the literature). In practice, both methods might converge
faster, and it remains open to study such faster convergence rates 
under additional suitable assumptions.

\section*{Acknowledgment}
The work was done during Da Tang's internship at Baidu Big  Data Lab
in Beijing. 
Tong Zhang would like to acknowledge NSF IIS-1250985, NSF IIS-1407939, 
and NIH R01AI116744 for supporting his research.
The authors would also like to thank Wotao Yin for helpful discussions.

\bibliographystyle{plain}
\bibliography{admm_note}

\appendix
\section{Proof of Theorem~\ref{thm:admm}}
The fact that $w^t$ minimizes the objective function in line 4 of
Algorithm~\ref{algo:admm}, together with the relationship of
$\alpha^t$ and $\alpha^{t-1}$ in line 5, implies that
\begin{equation}
\nabla \phi(w^{t}) + A^\top \alpha^{t} = H(w^{t-1}-w^t) . \label{eq:nabla-phi-admm}
\end{equation}
We thus obtain
\begin{align}
& \phi(w^{t}) - \phi(w) + (\alpha^{t})^\top A(w^{t}-w) 
 + (A w- Bv-c)^{\top} \tilde{H} A (w^{t-1}-w^{t}) \nonumber\\
\leq & - \frac{\gamma}{2} \|\nabla \phi(w^{t})-\nabla \phi(w)\|_2^2 \nonumber\\
& + \nabla \phi(w^{t})^\top (w^{t} - w) + (\alpha^{t})^\top A (w^{t}-w)  
 + (A w- Bv-c)^{\top} \tilde{H} A (w^{t-1}-w^{t}) \nonumber\\
=&  - \frac{\gamma}{2} \|A^\top(\alpha^t-\alpha) +
H(w^t-w^{t-1})\|_2^2 \nonumber\\
& + (w^{t-1}-w^{t})^\top A^\top \tilde{H} A (w^t-w)   + (A w- Bv-c)^{\top} \tilde{H} A (w^{t-1}-w^{t}) \nonumber\\
=&  - \frac{\gamma}{2} \|A^\top(\alpha^t-\alpha) +
H(w^t-w^{t-1})\|_2^2\nonumber\\
& + \frac{1}{2} \left[ \|Aw^{t-1}-Bv-c\|_{\tilde{H}}^2 -
  \|Aw^{t}-Bv-c\|_{\tilde{H}}^2 -
  \|Aw^t-Aw^{t-1}\|_{\tilde{H}}^2\right] .
\label{eq:dphi-admm}
\end{align}
In the above derivation, the
inequality is a direct consequence of the smoothness of $\phi$, which
implies that for any $w'$ and $w$, $\phi(w) \geq \phi(w') + \nabla
\phi(w')^\top (w-w') + 0.5 \gamma \|\nabla \phi(w)-\nabla
\phi(w')\|_2^2$. 
The first equality is due to \eqref{eq:nabla-phi-admm}, 
and $\nabla \phi(w) + A^\top \alpha = 0$ (which follows from
the assumption $w=\nabla \phi^*(-A^\top \alpha)$ of the theorem). 
The second equality is algebra.

We also have from the optimality of $v^t$ for minimizing the objective
function in line 3 of Algorithm~\ref{algo:admm}, and the relationship of
$\alpha^t$ and $\alpha^{t-1}$ in line 5:
\begin{equation}
\nabla g(v^{t}) - B^\top \alpha^t = -\rho G(v^t-v^{t-1}) + \rho B^\top A ({w}^{t-1}-w^t) . \label{eq:nabla-g-admm}
\end{equation}
Therefore
\begin{align}
&g(v^t) - g(v) + \frac{\lambda}{2} \|v^t-v\|_2^2- \alpha^{t\top} B (v^t - v) \nonumber\\
\leq& \nabla g(v^t)^\top (v^t-v) - \alpha^{t\top} B (v^t - v) \nonumber\\
=& \rho (v^t-v^{t-1})^\top G (v-v^t) + \rho (w^t-w^{t-1})^\top A^\top B(v-v^t) \nonumber \\
=& \frac{\rho}{2} [\|v-v^{t-1}\|_G^2 - \|v-v^{t}\|_G - \|v^t-v^{t-1}\|_G^2] \nonumber\\
& +
\frac{\rho}{2} [ \|A w^{t-1}-B v-c\|_2^2 + \|A w^t-Bv^t-c\|_2^2 -
\|Aw^{t}-Bv-c\|_2^2-\|Aw^{t-1}-Bv^{t}-c\|_2^2 ] \nonumber \\
=& \frac{\rho}{2} [\|v-v^{t-1}\|_G^2 - \|v-v^{t}\|_G -
\|v^t-v^{t-1}\|_G^2] 
+\frac{1}{2\rho} \|\alpha^t-\alpha^{t-1}\|_2^2
\nonumber \\
& +
\frac{\rho}{2} [ \|A w^{t-1}-B v-c\|_2^2 -
\|Aw^{t}-Bv-c\|_2^2-\|Aw^{t-1}-Bv^{t}-c\|_2^2 ] \label{eq:dg-admm} .
\end{align}
In the above derivation, the first inequality is due to the strong
convexity of $g(\cdot)$. The first equality employs
\eqref{eq:nabla-g-admm}.
The second equality is algebra, and the third equality is due to 
the relationship of $\alpha^t$ and $\alpha^{t-1}$ in line 5 of Algorithm~\ref{algo:admm}.

Finally we have
\begin{align}
&-(\alpha^t-\alpha)^\top (Aw^t-Bv^t-c)\nonumber\\
=&- \frac{1}{\rho} (\alpha^t-\alpha)^\top(\alpha^t-\alpha^{t-1})\nonumber\\
=& \frac{1}{2\rho}[
\|\alpha-\alpha^{t-1}\|_2^2-\|\alpha-\alpha^t\|_2^2-\|\alpha^t-\alpha^{t-1}\|_2^2] ,
\label{eq:dalpha-admm}
\end{align}
where the first equality uses 
the relationship of $\alpha^t$ and $\alpha^{t-1}$ in line 5 of
Algorithm~\ref{algo:admm}, and the second equality is algebra. 

By adding \eqref{eq:dphi-admm}, \eqref{eq:dg-admm}, \eqref{eq:dalpha-admm}, we obtain
\begin{align*}
&\phi(w^{t}) + g(v^t) - \phi(w) - g(v)
- \tilde{\alpha}^{t\top}(Aw-Bv-c)
+ \alpha^\top (Aw^{t}-Bv^t-c)\\
\leq&
  - \frac{\gamma}{2} \|A^\top(\alpha^t-\alpha) +H(w^t-w^{t-1})\|_2^2 
- \frac{\lambda}{2} \|v^t-v\|_2^2
\nonumber\\
& + \frac{1}{2} \left[ \|Aw^{t-1}-Bv-c\|_{\tilde{H}}^2 -
  \|Aw^{t}-Bv-c\|_{\tilde{H}}^2 -
  \|Aw^t-Aw^{t-1}\|_{\tilde{H}}^2\right] \nonumber\\
& + \frac{\rho}{2} [\|v-v^{t-1}\|_G^2 - \|v-v^{t}\|_G -
\|v^t-v^{t-1}\|_G^2] 
\nonumber \\
& +
\frac{\rho}{2} [ \|A w^{t-1}-B v-c\|_2^2 -
\|Aw^{t}-Bv-c\|_2^2-\|Aw^{t-1}-Bv^{t}-c\|_2^2 ] \nonumber\\
&+ \frac{1}{2\rho}[
\|\alpha-\alpha^{t-1}\|_2^2-\|\alpha-\alpha^t\|_2^2] ,
\end{align*}
which can be rewritten as the following bound:
\begin{align*}
r_t
\leq&
  \underbrace{- \frac{\gamma}{2} \|A^\top(\alpha^t-\alpha) +H(w^t-w^{t-1})\|_2^2 
- \frac12 \|w^t-w^{t-1}\|_{H}^2+ \frac{\theta}{2\rho}\|\alpha^t-\alpha\|_2^2}_{X_t}\\
&\underbrace{- \frac{\lambda}{2} \|v^t-v\|_2^2 + \frac{\rho\theta}{2}\|v-v^{t}\|_G}_{Y_t}
-\frac{\rho}{2}\|v^t-v^{t-1}\|_G^2 
\nonumber\\
& \underbrace{+ \frac{\rho\theta}{2(1+\theta)}\|A w^{t-1}-B v-c\|_2^2 
+ \frac{\theta}{2(1+\theta)} \|Aw^{t-1}-Bv-c\|_{\tilde{H}}^2 
-\frac{\rho}{2}\|Aw^{t-1}-Bv^{t}-c\|_2^2}_{Z_t} 
\nonumber\\
& + \frac{1}{2} \left[ \frac{1}{1+\theta} \|Aw^{t-1}-Bv-c\|_{\tilde{H}}^2 -
  \|Aw^{t}-Bv-c\|_{\tilde{H}}^2\right] \nonumber\\
& + \frac{\rho}{2} [\|v-v^{t-1}\|_G^2 - (1+\theta) \|v-v^{t}\|_G] 
\nonumber \\
& +
\frac{\rho}{2} \left[ \frac{1}{1+\theta} \|A w^{t-1}-B v-c\|_2^2 -
\|Aw^{t}-Bv-c\|_2^2\right] \nonumber\\
&+ \frac{1}{2\rho}[
\|\alpha-\alpha^{t-1}\|_2^2-(1+\theta)\|\alpha-\alpha^t\|_2^2 ] \\
=& X_t + Y_t -\frac{\rho}{2}\|v^t-v^{t-1}\|_G^2+ Z_t + (1+\theta)^{-1} \delta_{t-1} - \delta_t .
\end{align*}
We can bound $X_t$ as follows:
\begin{align*}
X_t=&- \frac{\gamma}{2} \|A^\top(\alpha^t-\alpha) +H(w^t-w^{t-1})\|_2^2 
-  \frac12 \|w^t-w^{t-1}\|_{H}^2+
\frac{\theta}{2\rho}\|\alpha^t-\alpha\|_2^2\\
\leq&- \frac{\gamma}{2} \|A^\top(\alpha^t-\alpha) +H(w^t-w^{t-1})\|_2^2 
-  \frac{1}{2\sMH} \|H(w^t-w^{t-1})\|_{2}^2+
\frac{\theta}{2\rho}\|\alpha^t-\alpha\|_2^2\\
\leq& \frac12 \max_u \left[ - \gamma
  \|A^\top(\alpha^t-\alpha) +u\|_2^2 -\sMH^{-1} \|u\|_2^2 \right]
+ \frac{\theta}{2\rho}\|\alpha^t-\alpha\|_2^2\\
=& - \frac{\gamma/2}{\gamma \sMH+1}  \|A^\top(\alpha^t-\alpha) \|_2^2 
+ \frac{\theta}{2\rho}\|\alpha^t-\alpha\|_2^2 \leq 0 .
\end{align*}
The last inequality uses the assumption on $\theta$ in the theorem. 
We also have
\begin{align*}
Z_t
=& {\frac{\theta}{2(1+\theta)}\|Aw^{t-1}- Bv-c\|_{\tilde{H}}^2 + 
\frac{\rho \theta}{2(1+\theta)} \|A w^{t-1}-B v-c\|_2^2
-\frac{\rho}{2} \|Aw^{t-1}-Bv^{t}-c\|_2^2} \\
\leq & \frac{\rho}{2}
\left[ {\frac{\theta (1+\stH/\rho)}{1+\theta}\|Aw^{t-1}- (Bv+c)\|_2^2 
- \|Aw^{t-1}-Bv^{t}-c\|_2^2} \right] \\
\leq& \frac{\theta \rho (\rho +\stH)}{2(\rho-\theta \stH)}\| B(v^t-v)\|_2^2 ,
\end{align*}
where the second inequality uses the fact that
\[
\frac{\theta (1+a)}{1+\theta} \|u\|_2^2 - \|u'\|_2^2 \leq
\frac{1+a}{1-\theta a} \theta \|u-u'\|_2^2 ,
\]
when $\theta a <1$ with $a=\stH/\rho$. Therefore
\begin{align*}
Y_t + Z_t \leq &
- \frac{\lambda}{2} \|v^t-v\|_2^2 + \frac{\rho\theta}{2}\|v-v^{t}\|_G
+ \frac{\theta \rho (\rho +\stH)}{2(\rho-\theta \stH)} \| B(v^t-v)\|_2^2 \\
\leq & 
\left[ - \frac{\lambda}{2} + \frac{\rho\theta}{2} \sG
+ \frac{\theta (\rho +\stH)}{2(1-s)} \sB^2 \right] \| v^t-v\|_2^2
\leq 0.
\end{align*}
Therefore we obtain 
\begin{align*}
r_t \leq X_t + Y_t + Z_t + (1+\theta)^{-1} \delta_{t-1} - \delta_t
\leq (1+\theta)^{-1} \delta_{t-1} - \delta_t .
\end{align*}
Now by multiplying the above displayed inequality by
$(1+\theta)^{t-T}$, and sum over $t=1,\ldots,T$, we obtain \eqref{eq:admm-primal}.

In order to obtain \eqref{eq:admm-dual}, we simply note that
\eqref{eq:nabla-phi-admm} implies that
\begin{equation}
\nabla \phi^*(-A^\top \tilde{\alpha}^{t}) - w^t = 0. \label{eq:nabla-phi-star-admm}
\end{equation}
Therefore \eqref{eq:dphi-admm} can be replaced by the following
inequality:
\begin{align}
& \phi^*(-A^\top \tilde{\alpha}^{t}) - \phi^*(-A^\top \alpha) +
(\alpha^{t}-\alpha)^\top A w^{t} 
 + (- Bv-c)^{\top} \tilde{H} A (w^{t-1}-w^{t}) \nonumber\\
\leq & \nabla \phi^*(-A^\top \tilde{\alpha}^{t})^\top (-A^\top\tilde{\alpha}^{t} + A^\top\alpha) + (\alpha^{t}-\alpha)^\top A 
w^{t} \nonumber\\
&- \frac{\gamma}{2} \|-A^\top \tilde{\alpha}^t + A^\top \alpha\|_2^2
 + (- Bv-c)^{\top} \tilde{H} A (w^{t-1}-w^{t}) \nonumber\\
=& - (w^t)^\top (A^\top {\alpha}^{t} + H(w^t-w^{t-1})- A^\top\alpha) + (\alpha^{t}-\alpha)^\top A 
w^{t} \nonumber\\
&- \frac{\gamma}{2} \|A^\top(\alpha^t-\alpha) +H(w^t-w^{t-1})\|_2^2
 + (- Bv-c)^{\top} \tilde{H} A (w^{t-1}-w^{t}) \nonumber\\
=&  - \frac{\gamma}{2} \|A^\top(\alpha^t-\alpha) +
H(w^t-w^{t-1})\|_2^2\nonumber\\
& + \frac{1}{2} \left[ \|Aw^{t-1}-Bv-c\|_{\tilde{H}}^2 -
  \|Aw^{t}-Bv-c\|_{\tilde{H}}^2 -
  \|Aw^t-Aw^{t-1}\|_{\tilde{H}}^2\right] ,
\label{eq:dphi-star-admm}
\end{align}
where the first inequality uses the fact $\phi^*$ is $\gamma$ strongly
convex, which is a direct consequence of the fact that $\phi$ is
$1/\gamma$ smooth. The first equality is due to
\eqref{eq:nabla-phi-star-admm} and the definition of
$\tilde{\alpha}^t$. The second equality is algebra.

Now, we note that the right hand side of
\eqref{eq:dphi-star-admm}
is the same as that of \eqref{eq:dphi-admm}.
Therefore the remaining of the proof follows the same argument as that of \eqref{eq:admm-primal},
where we simply use the addition of 
\eqref{eq:dphi-star-admm}, \eqref{eq:dg-admm}, and \eqref{eq:dalpha-admm}
to replace the addition of
\eqref{eq:dphi-admm}, \eqref{eq:dg-admm}, and \eqref{eq:dalpha-admm}.
This leads to \eqref{eq:admm-dual}.

\section{Proof of Corollary~\ref{cor:admm-duality-gap}}

We have from \eqref{eq:admm-dual}
\begin{align*}
&
2(1+\theta)^T \left[ \phi^*(-A^\top \tilde{\alpha}^{T}) + g(v^T) - \phi^*(-A^\top \alpha) - g(v)
+ (\tilde{\alpha}^{T})^{\top}(Bv+c)
- \alpha^\top (Bv^T+c)\right]\\
\leq&
(\rho+\stH) \|Aw^{0}-Bv-c\|_2^2 +
\rho(1+\theta) \|v^0-v\|_G^2 + \frac{1+\theta}{\rho}
\|\alpha-\alpha^{0}\|_2^2 \\
\leq&
2 (\rho+\stH) \|Aw^{0}-Bv^0-c\|_2^2 \\
& +
(2 (\rho+\stH) \sB^2+ \rho(1+\theta) \sG)\|v^0-v\|_2^2 
+ \frac{1+\theta}{\rho} \|\alpha-\alpha^{0}\|_2^2 .
\end{align*}

Now we set $\alpha= -(A^{+})^{\top} \nabla \phi (A^{+}(B v^T+c))$ and 
$v=\nabla g^*(B^\top \tilde{\alpha}^T)$.
This choice achieves the maximum
value of the
left hand side over $(\alpha,v)$. With this choice, and
the definition of convex conjugate, we obtain 
\begin{align}
&2(1+\theta)^T [\phi(A^{+}(B{v}^T+c)) +  g({v}^T) -
D(\tilde{\alpha}^T) ] \nonumber\\
\leq&
(2 (\rho+\stH) \sB^2+ \rho(1+\theta) \sG)\|v^0-v\|_2^2 
+ \frac{1+\theta}{\rho} \|\alpha-\alpha^{0}\|_2^2 \nonumber\\
&+ 2 (\rho+\stH) \|Aw^{0}-Bv^0-c\|_2^2 . \label{eq:dgap-estimate}
\end{align}

From Corollary~\ref{cor:admm-bregman}, we obtain
\begin{equation}
 \frac{\rho}{2} \|A(w^T-w_*)\|_2^2 
+ \frac{1+\theta}{2\rho} \|\alpha^T-\alpha_*\|_2^2
+ \frac{\rho(1+\theta)}{2} \|v^T-v_*\|_G^2
\leq   
\frac{(1+\theta)^{-T}}{2} \delta_*^0 . \label{eq:param-bound}
\end{equation}
Therefore 
\[
\|A(w^T-w^{T-1})\|_2^2\leq 2\|A(w^T-w)\|_2^2+2 \|A(w^{T-1}-w)\|_2^2\leq
2 (2+\theta) (1+\theta)^{-T}\delta_*^0/\rho .
\]
Moreover, \eqref{eq:param-bound} also implies $\|\alpha^T-\alpha_*\|_2^2\leq \rho(1+\theta)^{-1} (1+\theta)^{-T}\delta_*^0$.
Therefore 
\begin{align*}
&\|\tilde{\alpha}^T-\alpha_*\|_2 \leq \|\alpha^T-\alpha_*\|_2 + \stH
\|A(w^T-w^{T-1})\|_2\\
\leq& \stH \sqrt{2(2+\theta) (1+\theta)^{-T}\delta_*^0/\rho} +
\sqrt{\rho(1+\theta)^{-1} (1+\theta)^{-T}\delta_*^0} .
\end{align*}
It follows from the definition of $b_2(\cdot)$ that
\[
\|v-v^{0}\|_2^2\leq b_2((1+\theta)^{-T}\delta_*^0) .
\]
Similarly, we obtain from \eqref{eq:param-bound} that $\|v^T-v_*\|_G^2 \leq (1+\theta)^{-T}\delta_*^0/(\rho+\rho\theta)$. 
It implies that 
\[
\|\alpha-\alpha^{0}\|_2^2 \leq b_1((1+\theta)^{-T}\delta_*^0) .
\]
Now the first desired bound of the theorem can be obtained by plugging in the estimates
of
$\|v-v^{0}\|_2^2$ and $\|\alpha-\alpha^{0}\|_2^2$ into \eqref{eq:dgap-estimate}.

For the second desired bound, we note from the Jensen's inequality and
\eqref{eq:admm-dual} that
\begin{align*}
&
\left[ \phi^*(-A^\top \bar{\alpha}^{T}) + g(\bar{v}^T) - \phi^*(-A^\top \alpha) - g(v)
+ (\bar{\alpha}^{T})^{\top}(Bv+c)
- \alpha^\top (B\bar{v}^T+c)\right]\\
\leq&
\frac{1}{\sum_{t=1}^T (1+\theta)^{t-1}}
\left[
  \frac{1}{2} (\rho+\stH)\|Aw^{0}-Bv-c\|_2^2 +  \frac{\rho}{2}
  \|v^0-v\|_G^2 + \frac{1}{2\rho}
  \|\alpha-\alpha^{0}\|_2^2  
\right] . 
\end{align*}
Again we simply take the choice of $(\alpha,v)$ 
that achieves the maximum on the left hand side:
$\alpha= -(A^{+})^{\top} \nabla \phi (A^{+}(B \bar{v}^T+c))$ and 
$v=\nabla g^*(B^\top \bar{\alpha}^T)$.

\section{Proof of Theorem~\ref{thm:lin-admm}}

The basic proof structure  is the same as that of Theorem~\ref{thm:admm}. 
The fact that $w^t$ minimizes the objective function in line 4 of Algorithm~\ref{algo:lin-admm},
together with the relationship of $\alpha^t$ and $\alpha^{t-1}$ in line 5,
implies that
\begin{equation}
\nabla \phi(w^{t-1}) + A^\top \alpha^{t} = H(w^{t-1}-w^{t})  . \label{eq:beta-lin-admm}
\end{equation}
We thus obtain
\begin{align}
& \phi(w^{t-1}) - \phi(w) +(\alpha^{t})^\top A(w^{t}-w)  + \alpha^\top A(w^{t-1}-w^t) 
\nonumber\\
& + (A w- Bv-c)^{\top} \tilde{H} A (w^{t-1}-w^{t}) 
\nonumber\\
\leq & \nabla \phi(w^{t-1})^\top (w^{t-1} - w) 
+ (\alpha^{t})^\top A(w^{t}-w) + \alpha^\top A(w^{t-1}-w^t) \nonumber\\
& + (A w- Bv-c)^{\top} \tilde{H} A (w^{t-1}-w^{t}) 
-\frac{\gamma}{2} \|\nabla \phi(w^{t-1})-\nabla \phi(w)\|_2^2 \nonumber\\
=& (H(w^{t-1}-w^t)-A^\top \alpha^t)^\top (w^{t-1}-w)
+ (\alpha^{t})^\top A(w^{t}-w) + \alpha^\top A(w^{t-1}-w^t) \nonumber\\
& + (A w- Bv-c)^{\top} \tilde{H} A (w^{t-1}-w^{t}) 
 - \frac{\gamma}{2} \|A^\top(\alpha^t-\alpha) + H(w^t-w^{t-1})\|_2^2\nonumber\\
=&
(w^{t}-w^{t-1})^\top (A^\top(\alpha^t-\alpha) )
 - \frac{\gamma}{2} \|A^\top(\alpha^t-\alpha) + H(w^t-w^{t-1})\|_2^2\nonumber\\
& + \frac{1}{2} \left[ \|Aw^{t-1}-Bv-c\|_{\tilde{H}}^2 -
  \|Aw^{t}-Bv-c\|_{\tilde{H}}^2 +
  \|w^t-w^{t-1}\|_{{H}}^2\right] ,
\label{eq:dphi-lin-admm}
\end{align}
where the derivation uses similar arguments as those of
\eqref{eq:dphi-admm}.
The first inequality uses the smoothness of $\phi$, and the first
equality uses \eqref{eq:beta-lin-admm}. The second equality is algebra.

We also have from the optimality of $v^t$ for minimizing the objective
function in line 3 of Algorithm~\ref{algo:lin-admm}, and the relationship of
$\alpha^t$ and $\alpha^{t-1}$ in line 5, to obtain \eqref{eq:dg-admm}.
Finally, we can also obtain \eqref{eq:dalpha-admm}.

By adding \eqref{eq:dphi-lin-admm}, \eqref{eq:dg-admm},
\eqref{eq:dalpha-admm}, 
and use the simplified notation $\Delta w=(w^t-w^{t-1})$, and $\Delta \alpha=\alpha^t-\alpha$, 
we obtain 
\begin{align*}
r_t \leq&
\underbrace{\Delta w^\top (A^\top\Delta \alpha)
 - \frac{\gamma}{2} \|A^\top \Delta \alpha + H\Delta w\|_2^2
+ \frac12 \|\Delta w\|_{H}^2
+ \frac{\theta}{2\rho}\|\Delta \alpha\|_2^2}_{X_t}
\\
& \underbrace{- \frac{\lambda}{2} \|v^t-v\|_2^2 + 
\frac{\rho \theta}{2} \|v-v^{t}\|_G}_{Y_t} - \frac{\rho}{2}
\|v^t-v^{t-1}\|_G^2 \\
&
\underbrace{+\frac{\theta}{2(1+\theta)}\|w^{t-1}- A^{-1}(Bv+c)\|_H^2 + 
\frac{\rho \theta}{2(1+\theta)} \|A w^{t-1}-B v-c\|_2^2
-\frac{\rho}{2} \|Aw^{t-1}-Bv^{t}-c\|_2^2}_{Z_t} \\
&+
\frac{1}{2}
\left[\frac{1}{1+\theta}\|Aw^{t-1}- (Bv+c)\|_{\tilde{H}}^2 -\|Aw^t-(Bv+c)\|_{\tilde{H}}^2 \right]\\
&+ \frac{\rho}{2} [\|v-v^{t-1}\|_G^2 - (1+\theta)\|v-v^{t}\|_G ] \\
&+
\frac{\rho}{2} [ \frac{1}{1+\theta}\|A w^{t-1}-B v-c\|_2^2 -\|Aw^{t}-Bv-c\|_2^2 ] \nonumber \\
&+ \frac{1}{2\rho}[
\|\alpha-\alpha^{t-1}\|_2^2-(1+\theta) \|\alpha-\alpha^t\|_2^2] .
\end{align*}
We can bound $X_t$ as follows:
\begin{align*}
X_t =&-
(H\Delta w)^\top (\gamma I-H^{-1}) (A^\top\Delta \alpha)
 - \frac{\gamma}{2} (\|A^\top \Delta \alpha\|_2^2 + \|H\Delta w\|_2^2)
+ \frac12 \|\Delta w\|_{H}^2
+ \frac{\theta}{2\rho}\|\Delta \alpha\|_2^2\\
\leq&
(\gamma-1/\smH)\|H\Delta w\|_2 
\|A^\top\Delta \alpha\|_2
- \frac{\gamma-1/\smH}{2} \|H\Delta w\|_2^2
  - \frac{\gamma}{2} \|A^\top \Delta \alpha\|_2^2 
   + \frac{\theta}{2\rho}\|\Delta \alpha\|_2^2\\
\leq&
- \frac{1}{2\smH}\|A^\top \Delta \alpha\|_2^2
 + \frac{\theta}{2\rho}\|\Delta \alpha\|_2^2
 \leq 0 .
\end{align*}
The first inequality uses the assumption that $\gamma - \smH^{-1}\geq
0$ in the theorem, and norm inequalities.
The second inequality is obtained by taking the maximum over
$\|H\Delta w\|_2$.
The last inequality uses the assumptions on $\theta$. 
We also can use the same derivation as that of Theorem~\ref{thm:admm}
to show that $Y_t + Z_t \leq 0$. 
Therefore
\[
r_t \leq  X_t + Y_t -\frac{\rho}{2}\|v^t-v^{t-1}\|_G^2+ Z_t + (1+\theta)^{-1} \delta_{t-1} - \delta_t 
\leq (1+\theta)^{-1} \delta_{t-1} - \delta_t .
\]
We can multiply the above by $(1+\theta)^{t-T}$ and then sum over $t=1,\ldots$ to obtain
\eqref{eq:lin-admm-primal}.

Similarly we can prove a dual version of \eqref{eq:dphi-lin-admm}
below. 
The equation in \eqref{eq:beta-lin-admm} and the definition of
$\tilde{\alpha}^t$ in the theorem imply that 
\[
w^{t-1}= \nabla \phi^*(-A^\top \tilde{\alpha}^t) . 
\]
We thus have
\begin{align}
& \phi^*(-A^\top \tilde{\alpha}^{t}) - \phi^*(-A^\top \alpha) +
(\alpha^{t}-\alpha)^\top A w^{t} 
 + (- Bv-c)^{\top} \tilde{H} A (w^{t-1}-w^{t}) \nonumber\\
\leq& 
 \nabla \phi^*(-A^\top \tilde{\alpha}^{t})^\top
(-A^\top\tilde{\alpha}^{t} + A^\top\alpha) 
 - \frac{\gamma}{2} \|A^\top(\tilde{\alpha}^t-\alpha)\|_2^2 \nonumber\\
&+ 
(\alpha^{t}-\alpha)^\top A w^{t}
+ (- Bv-c)^{\top} \tilde{H} A (w^{t-1}-w^{t}) \nonumber\\
=& (w^{t-1})^\top (-(A^\top \alpha^t+H(w^t-w^{t-1}))+A^\top \alpha)
+(\alpha^{t}-\alpha)^\top A w^{t}\nonumber\\
& - \frac{\gamma}{2} \|A^\top(\alpha^t-\alpha) +H(w^t-w^{t-1})\|_2^2
+ (- Bv-c)^{\top} \tilde{H} A (w^{t-1}-w^{t}) \nonumber\\
=& 
(w^{t}-w^{t-1})^\top (A^\top(\alpha^t-\alpha))
 - \frac{\gamma}{2} \|A^\top(\alpha^t-\alpha) +H(w^t-w^{t-1})\|_2^2\nonumber\\
& + \frac{1}{2} \left[ \|Aw^{t-1}-Bv-c\|_{\tilde{H}}^2 -
  \|Aw^{t}-Bv-c\|_{\tilde{H}}^2 +
  \|w^t-w^{t-1}\|_{{H}}^2\right] .
\label{eq:dphi-star-lin-admm}
\end{align}
In the above derivation, the first inequality uses the strong
convexity of $\phi^*$, which follows from the smoothness of $\phi$.
The first equality uses the relationship of
 $\nabla \phi^*(-A^\top\tilde{\alpha}^t)$ and $w^{t-1}$ and the relationship of $\tilde{\alpha}^t$ and $\alpha^t$.
The last equality uses algebra.
 Note that the right hand side of 
\eqref{eq:dphi-lin-admm} and that of \eqref{eq:dphi-star-lin-admm} are
the same.
Therefore by adding \eqref{eq:dphi-star-lin-admm}, \eqref{eq:dg-admm},
\eqref{eq:dalpha-admm}, we obtain 
\[
r_t^* \leq X_t + Y_t -\frac{\rho}{2}\|v^t-v^{t-1}\|_G^2 + Z_t + (1+\theta)^{-1} \delta_{t-1} - \delta_t 
\leq (1+\theta)^{-1} \delta_{t-1} - \delta_t .
\]
We can multiply $(1+\theta)^{t-T}$
to both sides, and then sum over $t=1,\ldots$ to obtain \eqref{eq:lin-admm-dual}.

\end{document}